\documentclass[12pt]{amsart}
\usepackage{fullpage}
\newtheorem{theorem}{Theorem}
\newtheorem{lemma}{Lemma}
\newtheorem{cor}{Corollary}

\theoremstyle{definition}

\def\SBIMSMark#1#2#3{
 \font\SBF=cmss10 at 10 true pt
 \font\SBI=cmssi10 at 10 true pt
 \setbox0=\hbox{\SBF Stony Brook IMS Preprint \##1}
 \setbox2=\hbox to \wd0{\hfil \SBI #2}
 \setbox4=\hbox to \wd0{\hfil \SBI #3}
 \setbox6=\hbox to \wd0{\hss
             \vbox{\hsize=\wd0 \parskip=0pt \baselineskip=10 true pt
                   \copy0 \break%
                   \copy2 \break%
                   \copy4 \break}}
 \dimen0=\ht6   \advance\dimen0 by \vsize \advance\dimen0 by 8 true pt
                \advance\dimen0 by -\pagetotal
 \dimen2=\hsize \advance\dimen2 by .25 true in
  \setbox0=\hbox to 3.1 true in{
                \vbox to \ht6{\hsize=3 true in \parskip=0pt  \noindent  
         {\it Bull. London Math. Soc.}~{\bf 30} (1998) 62--66
                \vfill}}
  \ht0=0pt \dp0=0pt
 \ht6=0pt \dp6=0pt
 \setbox8=\vbox to \dimen0{\vfill \hbox to \dimen2{\copy0 \hss \copy6}}
 \ht8=0pt \dp8=0pt \wd8=0pt
 \copy8
 \message{*** Stony Brook IMS Preprint #1, #2 ***}
}

\begin{document}
\title{On Critical Points of Proper Holomorphic Maps on the Unit disk}
\author{Saeed Zakeri}
\address{Department of Mathematics, SUNY at Stony Brook, NY 11794}
\email{zakeri@math.sunysb.edu}
\begin{abstract}
We prove that a proper holomorphic map on the unit disk in the complex plane
is uniquely determined up to post-composition with a M\"{o}bius transformation
by its critical points.
\end{abstract}
\maketitle

\SBIMSMark{1996/7}{June 1996}{}
\thispagestyle{empty}

This note will give a brief proof of the following known theorem:
\begin{theorem}
Let $c_1, \cdots, c_d$ be given (not necessarily distinct) points in the open
unit disk {\bf D} in the complex plane. Then there exists a unique proper
holomorphic map $f:{\bf D}\rightarrow {\bf D}$ of degree $d+1$ normalized as $f(0)=0$ and $f(1)=1$ with critical points at the $c_j$.
\end{theorem}
It is to be understood that whenever some $c_j$ repeats $m$ times, the corresponding map has local degree $m+1$ at $c_j$. 

Since given any two points $a$ and $b$ with $|a|<1$ and $|b|=1$ there exists
a unique conformal automorphism of the unit disk which maps $a$ to $0$ and $b$ to 1, we have the following version of the uniqueness part of the above theorem:
\begin{cor}
Two proper holomorphic maps $f,g:{\bf D}\rightarrow {\bf D}$ have the same critical points, counted with multiplicity, if and only if $f=\tau \circ g$ for
some conformal automorphism $\tau$ of the unit disk. 
\end{cor} 
From the point of view of complex analysis, it is quite natural to ask questions about the dependence of such maps on their critical points, but this kind of question is also of some interest when one studies the parameter space for complex polynomial maps of a given degree. For example in \cite{M}, this theorem
immediately implies that the space of all ``critically marked'' normalized Blaschke products of degree $d+1$ is a topological cell of real dimension $2d$.
Since these provide a model space for hyperbolic polynomial maps with a given post-critical pattern, the discussion of hyperbolic components with marked critical points would become much easier and more natural having known this theorem.

The corresponding questions for the two other simply-connected Riemann surfaces (i.e., the complex plane and the Riemann sphere) have trivial and surprising 
answers. Every proper holomorphic map on the plane is a polynomial, so given a finite number of points in the plane one can always find a polynomial with corresponding critical points, and this polynomial is unique up to post-composition with a complex affine transformation. In the case of the Riemann sphere, however, both existence and uniqueness parts of the theorem are false. Every proper holomorphic map of degree $d+1$ on the sphere is a rational
map with 2$d$ critical points, counted with multiplicity. So, for example, it
is impossible to realize a single point with multiplicity~2 as the critical
set  of a rational map, since any such map would have degree~2 and local
degree~3 
near the double critical point. However, it can be shown that any $2d$ {\it distinct} points on the sphere can be realized as the critical set of a degree
$d+1$ rational map, and typically there is more than one way to get such a map,
even up to post-composition with a M\"{o}bius transformation \cite{G}.

The possibility of such a theorem was first suggested to us by J. Milnor. Our efforts to locate references led to a Chinese language paper \cite{W-P} in which the result is stated as a corollary of several lemmas without proofs. 
There one can also find a reference to an older paper by M. Heins \cite{H} in which
he actually proves the Theorem among other things by using a general machinery
he develops to study conformal metrics on Riemann surfaces. Since the Theorem
is important in its own way, we believe it deserves to be proved directly using
a clear argument. After all, it is rather surprising that such a nice fact,
which has to be classical, is not well-known. After this paper was written,
I was informed that the same problem is discussed in T. Bousch's thesis (Sur quelques probl\'{e}ms de dynamique holomorphe, University Paris 11, Orsay,
1992). There he transfers the problem to the upper-half plane and uses purely
algebraic methods to show that the mapping $\Phi$ in this paper is a local 
diffeomorphism.\\

It follows from elementary complex analysis that every proper holomorphic map $f:{\bf D} \rightarrow
{\bf D}$ of degree $d+1$ is a finite {\it Blaschke product} of the form
\begin{equation}
f(z)=\lambda \prod_{j=1}^{d+1} \left( \frac{z-a_j}{1-\bar{a_j}z} \right )
\end{equation}
where $|a_j|<1$ and $|\lambda|=1$. Therefore such a map is uniquely determined
up to a rotation by its zeros $\{a_1, \cdots ,a_{d+1} \}$. Note that $f$ extends holomorphically to a neighborhood of the closed unit disk $\bar {\bf D}$, mapping the boundary circle $ |z|=1$ to itself. We can always normalize $f$ so that $f(0)=0$ (say $a_{d+1}=0$) and $f(1)=1$ by post-composing $f$ with appropriate conformal automorphisms of the disk. Therefore every normalized map in this sense can be written as
\begin{equation}
f(z)=z\prod_{j=1}^d \left ( \frac{1-\bar {a_j}}{1-a_j} \right ) \left ( \frac{z-a_j}{1-\bar{a_j}z} \right ),
\end{equation}
where $|a_j|<1$. Note that the $j$-th factor
\begin{equation} 
\beta_{a_j}(z)=\left ( \frac{1-\bar{a_j}}{1-a_j} \right ) \left ( \frac {z-a_j}{1-\bar{a_j}z} \right )
\end{equation} 
is the unique 
conformal automorphism of the disk which maps $a_j$ to 0 and fixes 1. It is easy to see that {\it if $\{ a_n\} $ is a sequence of points in {\bf D} which converges to
a boundary point $a$ with $|a|=1$, then the corresponding sequence of automorphisms $\{ \beta_{a_n}\} $ has a subsequence which converges to 
a constant function whose absolute value is 1, and the convergence is uniform
on compact subsets of the open unit disk.} In fact, if $a\neq 1$, it is easy to see that the whole sequence actually converges to the constant function 1. On the other hand, if $a=1$, it follows from (3) that for every accumulation point $\beta $ of the sequence $\{ \frac{1-\bar{a}_n}{1-a_n}\} $ on the unit circle one can find a subsequence of $\{ \beta_{a_n}\} $ which converges locally uniformly to $-\beta $. \\

We denote the space of all normalized Blaschke products of the form (2) by $\mathcal{B}_d$.\\

Let $\Sigma_d$ be the quotient of the $d$-fold product ${\bf D}^d$ by the action
of the symmetric group $S_d$ on $d$ letters. This is the space of all unordered $d$-tuples $\{a_1, \cdots, a_d \}$ of points in the unit disk. Note that there is a bijection $\mathcal{B}_d \rightarrow \Sigma_d$ mapping a normalized Blaschke product to the (unordered) zero set $\{a_1, \cdots, a_d \}$. The topology of $\Sigma_d$ coming from ${\bf D}^d$ will match the topology of $\mathcal{B}_d$ (uniform convergence on compact subsets of the disk) under this bijection, so that the identification is in fact a homeomorphism.

Now consider the mapping
$$\Phi :\Sigma_d \rightarrow \Sigma_d$$
which sends a normalized Blaschke product with zero set ${\bf a}=\{ a_1, \cdots, a_d \}$ to its (unordered) critical set ${\bf c}=\{ c_1, \cdots, c_d\}$. We shall prove that $\Phi$ is a homeomorphism. The Theorem will follow 
immediately.

We list the properties of $\Phi$ in the following lemmas. First of all, it is
clear that
\begin{lemma} 
$\Phi$ is continuous.
\end{lemma}

A simple compactness argument proves the following 
\begin{lemma}
$\Phi$ is a proper map.
\end{lemma}
\begin{proof}
Let the sequence ${\bf a}^n= \{a_1^n, \cdots, a_d^n\}$ in $\Sigma_d$ leave every compact subset of $\Sigma_d$. After passing to a subsequence and relabeling, we can assume that ${\bf a}^n$ converges to
some ${\bf a}= \{ a_1, \cdots, a_d\}$, where $|a_j|\leq 1$ and $|a_1|=1$.
It follows from the discussion after (3) that the corresponding sequence of normalized Blaschke products $f_n$ has a subsequence which converges locally uniformly to a finite Blaschke product $f$
of {\it lower} degree $d'+1<d+1$, with zeros at 0 and those $a_j$ with $|a_j|<1$. Since $f$ has $d'$ critical points in the unit disk counted
with multiplicity, it follows that the sequence $\Phi ({\bf a}^n)$ has to leave every compact subset of $\Sigma_d$, too. 
\end{proof}

Now let $f$ be any proper holomorphic map of the form (1) and consider the density $\sigma_f$ of the pull-back of the Poincar\'{e} metric $\frac {2|dz|}{1-|z|^2}$ on {\bf D} under $f$:\\
$$\sigma_f(z)=\frac{2|f'(z)|}{1-|f(z)|^2}.$$
The metric $\sigma_f(z)|dz|$ is conformal of constant curvature $-1$ away from
the critical points of $f$. The condition on the curvature takes the form
\begin {equation}
\Delta \log \sigma_f(z)=\sigma_f^2(z).
\end{equation}
By the Ahlfors-Schwarz lemma, $f$ decreases the Poincar\'{e} distances, so $\frac{(1-|z|^2)|f'(z)|}{1-|f(z)|^2}\leq 1$. This suggests that we consider the {\it distance-ratio} function
\begin{equation}
 R_f(z)= \frac{1}{2}(1-|z|^2)\sigma_f(z)=\frac{1-|z|^2}{1-|f(z)|^2}|f'(z)|
\end{equation}
which compares the pull-back metric with the Poincar\'{e} metric on the disk. 
It is easy to prove the following lemma which lists the basic properties of $R_f$:
\begin{lemma}
Let $R_f$ be defined as in (5). Then
\begin{enumerate}
\item $R_f\leq 1$, with equality if and only if $f$ is a conformal automorphism of the disk.
\item $R_{(f\circ g)}=(R_f\circ g)R_g$. In particular, for every conformal automorphism $\tau$ of the disk, $R_{\tau \circ f}=R_f$ and $R_{f \circ \tau}=R_f \circ \tau$.
\item $R_f$ is a non-negative function on the disk with zeros at the critical points of $f$. It is real-analytic away from these critical points. Moreover, if $c$ is a critical point of $f$ and if $f$ has local degree $m\geq 2$ at $c$, then $R_f(z)=|z-c|^{m-1} \tilde{R}(z)$, where $\tilde{R}(z)$ is 
real-analytic and positive in a neighborhood of $c$. 
\end{enumerate}
\end{lemma}
\begin{lemma}
$R_f$ has a continuous extension to the closed unit disk.
In fact, $\lim_{|z|\rightarrow 1}R_f(z)=1$.
\end{lemma}
\begin{proof}
We carry out the proof when $z\rightarrow 1$. The general case follows from  Lemma 3(2). First note that for any Blaschke product $f$ of the form (1), the logarithmic derivative $d(\log f)/d(\log z)= zf'/f$ is a rational function which is positive and equal to $|f'|$ on the unit circle $|z|=1$. In fact, direct computation shows that if $|z|=1$, then
$$\frac{zf'(z)}{f(z)}=\sum_{j=1}^{d+1}\frac{(1-|a_j|^2)z}{(z-a_j)(1-\bar{a_j}z)}= \sum_{j=1}^{d+1}\frac{1-|a_j|^2}{|z-a_j|^2}.$$
We proceed by induction on the degree of $f$. For $d=0$, $f$ is a conformal 
automorphism of the disk and the result is clear by Lemma 3(1). Suppose that the Lemma is true for all Blaschke products of degree $<d+1$.
Given a Blaschke product $f$ of degree $d+1$, use (1) to write $f=AB$, where $A$ and $B$ are Blaschke products of degree $<d+1$. Write
$$\frac{1-|f(z)|^2}{1-|z|^2}=\frac{(1-|A(z)|^2)|B(z)|^2+(1-|B(z)|^2)}{1-|z|^2}$$
By the induction hypothesis, the right side approaches $|A'(1)|+|B'(1)|$ as $z\rightarrow 1$. But $A$ and $B$ have positive logarithmic derivatives on the unit circle, so 
$$|A'(1)|+|B'(1)|=\frac{d(\log A)}{d(\log z)}(1)+\frac{d(\log B)}{d(\log z)}(1)=\frac{d(\log f)}{d(\log z)}(1)=|f'(1)|.$$ 
\end{proof}

{\it Remark.} Here is another proof for Lemma 4, as suggested to us by M. Lyubich: Consider the round annulus $B_r=\{ z: r<|z|<1\}$. For $r>0$ sufficiently close to $1$, $f$ is a $(d+1)$-to-1 covering map from $A_r=f^{-1}(B_r)$ to
$B_r$, and $A_r$ is a topological annulus with the unit circle as the outer 
boundary. Since $f$ is a covering, it is an isometry in the Poincar\'{e} metrics on $A_r$ and $B_r$. But a simple calculation shows that these Poincar\'{e} metrics are both asymptotic to the Poincar\'{e} metric on the unit disk as one approaches the unit circle. Hence the distance-ratio function has to tend to 1 near the unit circle.

\goodbreak
\begin{lemma} 
$\Phi$ is one-to-one.
\end{lemma}
\begin{proof}
Let two points in $\Sigma_d$ have the same image under $\Phi$. This means that the corresponding normalized Blaschke products $f$ and $g$ in $\mathcal{B}_d$ have the same critical points with the same multiplicities. Consider the two functions $\sigma_f$ and $\sigma_g$ and form the function 
$$h(z)=\frac{\sigma_f(z)}{\sigma_g(z)}=\frac{R_f(z)}{R_g(z)}$$
which is real-analytic and positive away from the critical points of $f$ (or $g$ ). But by (5) and Lemma 3(3), it is easy to see that $h$ can be extended to a real-analytic positive function on the whole unit disk. Also by Lemma 4, $h(z)\rightarrow 1$ as $|z|\rightarrow 1$. Let us show that $h\leq 1$. Suppose,
on the contrary, that $h(z)>1$ for some $z\in {\bf D}$. Then there is a point
$p\in {\bf D}$ at which $h$ has a maximum $h(p)>1$. Therefore $\Delta \log h(p)\leq 0$. On the other hand, by (4), $\Delta \log h(z)= \Delta \log \sigma_f(z) -\Delta \log \sigma_g(z)=\sigma_f^2(z)-\sigma_g^2(z)$ away from the critical points of $f$ (or $g$), hence everywhere in the unit disk by continuity. So $\sigma_f^2(p)-\sigma_g^2(p)\leq 0$, or $\sigma_f(p)\leq \sigma_g(p)$, which is a contradiction. Therefore $h\leq 1$.\\
A similar argument shows that $1/h\leq 1$. So $h=1$ everywhere in the disk. This
implies 
$$\frac{|f'(z)|}{1-|f(z)|^2}=\frac{|g'(z)|}{1-|g(z)|^2}.$$
But this means that $g(z)\mapsto f(z)$ is a Poincar\'{e} isometry as $z$ varies over any small region without critical points. Extending this isometry to a
conformal automorphism $\tau$ of the disk, we have $f=\tau \circ g$ everywhere
by analytic continuation. Since $f$ and $g$ fix 0 and 1, $\tau$ is the identity
and we have $f=g$. 
\end{proof}

\begin{cor}
$\Phi$ is a homeomorphism.
\end{cor}

\begin{proof}
$\Phi$ is a one-to-one continuous proper map $\Sigma_d \rightarrow
\Sigma_d$. Hence it is a covering map of degree 1.
\end{proof}

{\it Acknowledgement.} I would like to thank J. Milnor and M. Lyubich for their useful suggestions.

\bibliographystyle{amsplain}

\end{document}